\documentclass[envcountsame]{llncs}

\usepackage{amsmath}
\usepackage{amssymb}
\usepackage{pstricks}

\newcommand{\C}{\mathbb C}
\newcommand{\F}{\mathbb F}

\newcommand{\PP}{\mathbb P}
\newcommand{\Q}{\mathbb Q}
\newcommand{\R}{\mathbb R}
\newcommand{\Z}{\mathbb Z}

\newcommand{\frakb}{\mathfrak b}
\newcommand{\frakc}{\mathfrak c}

\newcommand{\frakD}{\mathfrak D}

\newcommand{\frakm}{\mathfrak m}
\newcommand{\frakN}{\mathfrak N}
\newcommand{\frakp}{\mathfrak p}

\newcommand{\calH}{\mathcal H}
\newcommand{\calO}{\mathcal O}

\newcommand{\alphahat}{\widehat{\alpha}}
\newcommand{\bhat}{\widehat{b}}
\newcommand{\betahat}{\widehat{\beta}}
\newcommand{\Bhat}{\widehat{B}}
\newcommand{\Fhat}{\widehat{F}}
\newcommand{\calOhat}{\widehat{\mathcal O}}
\newcommand{\gammahat}{\widehat{\gamma}}

\newcommand{\deltahat}{\widehat{\delta}}
\newcommand{\pihat}{\widehat{\pi}}
\newcommand{\phat}{\widehat{p}}
\newcommand{\varpihat}{\widehat{\varpi}}

\newcommand{\Zhat}{\widehat{\mathbb Z}}
\newcommand{\ZFhat}{\widehat{\mathbb Z}_F}

\DeclareMathOperator{\Cl}{Cl}
\DeclareMathOperator{\Coind}{Coind}
\DeclareMathOperator{\Gal}{Gal}
\DeclareMathOperator{\Hom}{Hom}
\DeclareMathOperator{\impart}{Im}
\DeclareMathOperator{\M}{M}
\DeclareMathOperator{\N}{N}
\DeclareMathOperator{\nrd}{nrd}
\DeclareMathOperator{\PGL}{PGL}

\DeclareMathOperator{\sgn}{sgn}

\newcommand{\slsh}[1]{\,|_{#1}\,}
\newcommand{\quat}[2]{\displaystyle{\biggl(\frac{#1}{#2}\biggr)}}

\begin{document}

\title{Computing automorphic forms on Shimura curves over fields with arbitrary class number}
\author{John Voight}
\institute{Department of Mathematics and Statistics \\ University of Vermont \\ 16 Colchester Ave \\ Burlington, VT\ 05401, USA \\ \email{jvoight@gmail.com}}

\maketitle

\begin{abstract}
We extend methods of Greenberg and the author to compute in the cohomology of a Shimura curve defined over a totally real field with arbitrary class number.  Via the Jacquet-Langlands correspondence, we thereby compute systems of Hecke eigenvalues associated to Hilbert modular forms of arbitrary level over a totally real field of odd degree.  We conclude with two examples which illustrate the effectiveness of our algorithms.
\end{abstract}

The development and implementation of algorithms to compute with automorphic forms has emerged as a major topic in explicit arithmetic geometry.  The first such computations were carried out for elliptic modular forms, and now very large and useful databases of such forms exist \cite{CremonaTables,Stein,SteinWatkins}.  Recently, effective algorithms to compute with Hilbert modular forms over a totally real field $F$ have been advanced.  The first such method is due to Demb\'el\'e \cite{Dembele0,Dembele}, who worked initially under the assumption that $F$ has even degree $n=[F:\Q]$ and strict class number $1$.  Exploiting the Jacquet-Langlands correspondence, systems of Hecke eigenvalues can be identified inside spaces of automorphic forms on $B^\times$, where $B$ is the quaternion algebra over $F$ ramified precisely at the infinite places of $F$---whence the assumption that $n$ is even.  Demb\'el\'e then provides a computationally efficient theory of Brandt matrices associated to $B$.  This method was later extended (in a nontrivial way) to fields $F$ of arbitrary class number by Demb\'el\'e and Donnelly \cite{DD}.  

When the degree $n$ is odd, a different algorithm has been proposed by Greenberg and the author \cite{GV}, again under the assumption that $F$ has strict class number $1$.  This method instead locates systems of Hecke eigenvalues in the (degree one) cohomology of a Shimura curve, now associated to the quaternion algebra $B$ ramified at all but one real place and no finite place.  This method uses in a critical way the computation of a fundamental domain and a reduction theory for the associated quaternionic unit group  \cite{V-fd}; see Section 1 for an overview.  In this article, we extend this method to the case where $F$ has arbitrary (strict) class number.  Our main result is as follows; we refer the reader to Sections 1 and 2 for precise definitions and notation.  

\begin{theorem} \label{maintheorem}
There exists an (explicit) algorithm which, given a totally real field $F$ of degree $n=[F:\Q]$, a quaternion algebra $B$ over $F$ ramified at all but one real place, an ideal $\frakN$ of $F$ coprime to the discriminant $\frakD$ of $B$, and a weight $k \in (2\Z_{>0})^n$, computes the system of eigenvalues for the Hecke operators $T_\frakp$ with $\frakp\nmid \frakD\frakN$ and the Atkin-Lehner involutions $W_{\frakp^e}$ with $\frakp^e\parallel \frakD\frakN$ acting on the space of quaternionic modular forms $S_k^B(\frakN)$ of weight $k$ and level $\frakN$ for $B$.
\end{theorem}

In other words, there exists an explicit finite procedure which takes as input the field $F$, its ring of integers $\Z_F$, a quaternion algebra $B$ over $F$, an ideal $\frakN \subset \Z_F$, and the vector $k$ encoded in bits (each in the usual way), and outputs a finite set of number fields $E_f \subset \overline{\Q}$ and sequences $(a_f(\frakp))_\frakp$ encoding the Hecke eigenvalues for each cusp form constituent $f$ in $S_k^B(\frakN)$, with $a_f(\frakp) \in E_f$.

From the Jacquet-Langlands correspondence, applying the above theorem to the special case where $\frakD=(1)$ (and hence $n=[F:\Q]$ is odd), we have the following corollary.

\begin{corollary}
There exists an algorithm which, given a totally real field $F$ of odd degree $n=[F:\Q]$, an ideal $\frakN$ of $F$, and a weight $k \in (2\Z_{>0})^n$, computes the system of eigenvalues for the Hecke operators $T_\frakp$ and Atkin-Lehner involutions $W_{\frakp^e}$ acting on the space of Hilbert modular cusp forms $S_k(\frakN)$ of weight $k$ and level $\frakN$.
\end{corollary}

This corollary is not stated in its strongest form: in fact, our methods overlap with the methods of Demb\'el\'e and his coauthors whenever there is a prime $\frakp$ which exactly divides the level; see Remark \ref{dembeleoverlap} for more detail.  Combining these methods, Donnelly and the author \cite{DembVoight} are systematically enumerating tables of Hilbert modular forms, and the details of these computations (including the dependence on the weight, level, and class number, as well as a comparison of the runtime complexity of the steps involved) will be reported there \cite{DembVoight}, after further careful optimization.

A third technique to compute with automorphic forms, including Hilbert modular forms, has been advanced by Gunnells and Yasaki \cite{GY}.  They instead use the theory of Vorono\u{\i} reduction and sharbly complexes; their work is independent of either of the above approaches.  

This article is organized as follows.  In Section 1, we give an overview of the basic algorithm of Greenberg and the author which works over fields $F$ with strict class number $1$.  In Section 2, using an adelic language we address the complications which arise over fields of arbitrary class number, and in Section 3 we make this theory concrete and provide the explicit algorithms announced in Theorem \ref{maintheorem}.  Finally, in Section 4, we consider two examples, one in detail; our computations are performed in the computer system \textsf{Magma} \cite{Magma}.

The author would like to thank Steve Donnelly and Matthew Greenberg for helpful discussions as well as the referees for their comments.  The author was supported by NSF Grant No.~DMS-0901971.

\section{An overview of the algorithm for strict class number $1$}

In this section, we introduce the basic algorithm of Greenberg and the author \cite{GV} with a view to extending its scope to base fields of arbitrary class number; for further reading, see the references contained therein.

Let $F$ be a totally real field of degree $n=[F:\Q]$ with ring of integers $\Z_F$.  Let $F_+^\times$ be the group of totally positive elements of $F$ and let $\Z_{F,+}^\times = \Z_F^\times \cap F_+^\times$.  Let $B$ be a quaternion algebra over $F$ of discriminant $\frakD$.  Suppose that $B$ is split at a unique real place $v_1$, corresponding to an embedding $\iota_\infty:B \hookrightarrow B \otimes \R \cong \M_2(\R)$, and ramified at the other real places $v_2,\dots,v_n$.  Let $\calO(1) \subset B$ be a maximal order and let 
\[ \calO(1)_+^\times=\{\gamma \in \calO(1)^\times : v_1(\nrd(\gamma))>0\}=\{ \gamma \in \calO(1) : \nrd(\gamma) \in \Z_{F,+}^\times\} \] denote the group of units of $\calO(1)$ with totally positive reduced norm.  Let 
\[ \Gamma(1) = \iota_\infty( \calO(1)_+^\times / \Z_F^\times) \subset \PGL_2(\R)^+, \] 
so that $\Gamma(1)$ acts on the upper half-plane $\calH = \{z \in \C : \impart(z) > 0\}$ by linear fractional transformations.  Let $\frakN \subset \Z_F$ be an ideal coprime to $\frakD$, let $\calO=\calO_0(\frakN)$ be an Eichler order of level $\frakN$, and let $\Gamma=\Gamma_0(\frakN) = \iota_\infty( \calO_0(\frakN)_+^\times/ \Z_F^\times)$.  

Let $k=(k_1,\ldots,k_n) \in (2\Z_{>0})^n$ be a weight vector; for example, the case $k=(2,\dots,2)$ of parallel weight $2$ is of significant interest.  Let $S_k^B(\frakN)$ denote the finite-dimensional $\C$-vector space of quaternionic modular forms of weight $k$ and level $\frakN$ for $B$.  Roughly speaking, a form $f \in S_k^B(\frakN)$ is an analytic function $f:\calH \to W_k(\C)$ which is invariant under the weight $k$ action by the group $\gamma\in \Gamma$, where $W_k(\C)$ is an explicit right $B^\times$-module \cite[(2.4)]{GV} and $W_k(\C)=\C$ when $k$ is parallel weight $2$.  The space $S_k^B(\frakN)$ comes equipped with the action of Hecke operators $T_\frakp$ for primes $\frakp \nmid \frakD\frakN$ and Atkin-Lehner involutions $W_{\frakp^e}$ for prime powers $\frakp^e \parallel \frakD\frakN$.

The Jacquet-Langlands correspondence \cite[Theorem 2.9]{GV} (see Hida \cite[Proposition 2.12]{HidaCM}) gives an isomorphism of Hecke modules
\[ S_k^B(\frakN) \xrightarrow{\sim} S_k(\frakD\frakN)^{\text{$\frakD$-new}}, \]
where $S_k(\frakD\frakN)^{\text{$\frakD$-new}}$ denotes the space of Hilbert modular cusp forms of weight $k$ and level $\frakD\frakN$ which are new at all primes dividing $\frakD$.  Therefore, as Hecke modules one can compute equivalently with Hilbert cusp forms or with quaternionic modular forms.

We compute with the Hecke module $S_k^B(\frakN)$ by identifying it as a subspace in the degree one cohomology of $\Gamma(1)$, as follows. 
Let $V_k(\C)$ be the subspace of the algebra $\C[x_1,y_1,\ldots,x_n,y_n]$ consisting of those polynomials $q$ which are homogeneous in $(x_i,y_i)$ of degree $w_i=k_i-2$.  Then $V_k(\C)$ has a right action of the group $B^\times$ given by
\begin{equation} \label{Vk}
q^\gamma(x_1,y_1,\ldots,x_n,y_n)=\left(\prod_{i=1}^n(\det\gamma_i)^{-w_i/2}\right)
q((x_1\,\,y_1)\overline{\gamma}_1,\ldots,(x_n\,\,y_n)\overline{\gamma}_n)
\end{equation}
for $\gamma \in B^\times$, where $\overline{\phantom{x}}$ denotes the standard involution (conjugation) on $B$ and $\gamma_i=v_i(\gamma) \in \M_2(\C)$.  By the theorem of Eichler and Shimura \cite[Theorem 3.8]{GV}, we have an isomorphism of Hecke modules
\[
S_k^B(\frakN) \xrightarrow{\sim} H^1\bigl(\Gamma,V_k(\C)\bigr)^+
\]
where the group cohomology $H^1$ denotes the (finite-dimensional) $\C$-vector space of crossed homomorphisms $f:\Gamma \to V_k(\C)$ modulo coboundaries and $+$ denotes the $+1$-eigenspace for complex conjugation.  By Shapiro's lemma \cite[\S 6]{GV}, we then have a further identification
\begin{equation} \label{VkC}
S_k^B(\frakN) \xrightarrow{\sim} H^1\bigl(\Gamma,V_k(\C)\bigr)^+ \cong H^1(\Gamma(1), V(\C))^{+},
\end{equation}
where $V(\C)=\Coind_{\Gamma}^{\Gamma(1)} V_k(\C)$.

In the isomorphism (\ref{VkC}), the Hecke operators act as follows.  Let $\frakp$ be a prime of $\Z_F$ with $\frakp \nmid \frakD\frakN$ and let $\F_\frakp$ denote the residue class field of $\frakp$.  Since $F$ has strict class number $1$, by strong approximation \cite[Theor\`eme III.4.3]{Vigneras} there exists $\pi \in \calO$ such that $\nrd \pi$ is a totally positive generator for $\frakp$.  It follows that there are elements $\gamma_a \in \calO_+^\times$, indexed by $a\in \PP^1(\F_\frakp)$, such that
\begin{equation} \label{heckequat}
\calO_+^\times\/\pi\/\calO_+^\times = \bigsqcup_{a\in\PP^1(\F_\frakp)} \calO_+^\times\alpha_a
\end{equation}
where $\alpha_a= \pi\gamma_a$.  

Let $f:\Gamma(1) \to V(\C)$ be a crossed homomorphism, and let $\gamma \in \Gamma(1)$.  The decomposition (\ref{heckequat}) extends to $\calO(1)$ as
\[ \calO(1)_+^\times\/\pi\/\calO(1)_+^\times = \bigsqcup_{a\in\PP^1(\F_\frakp)} \calO(1)_+^\times\alpha_a. \]
Thus, there are elements $\delta_a \in \calO(1)_+^\times$ for $a \in \PP^1(\F_\frakp)$ and a unique permutation $\gamma^*$ of $\PP^1(\F_\frakp)$ such that
\begin{equation} \label{deltaa}
\alpha_a\gamma = \delta_a\alpha_{\gamma^*a}
\end{equation}
for all $a$.  We then define $f \slsh{} T_\frakp:\Gamma(1) \to V(\C)$ by
\begin{equation}\label{E:heckeformula}
(f \slsh{} T_\frakp)(\gamma) = \sum_{a\in\PP^1(\F_\frakp)}f(\delta_a)^{\alpha_a}.
\end{equation}
The space $S_k^B(\frakN)$ similarly admits an action of Atkin-Lehner operators $W_{\frakp^e}$ for primes $\frakp^e \parallel \frakD\frakN$.  

From this description, we see that the Hecke module $H^1(\Gamma(1),V(\C))^+$ is amenable to explicit computation.  First, we compute a finite presentation for $\Gamma(1)$ with a minimal set of generators $G$ and a solution to the word problem for the computed presentation using an algorithm of the author \cite{V-fd}.  Given such a set of generators and relations, one can explicitly find a basis for the $\C$-vector space $H^1(\Gamma(1),V(\C))$ \cite[\S 5]{GV}.  

We then compute the action of the Hecke operator $T_\frakp$ on $H^1(\Gamma(1),V(\C))$.  We first compute a splitting $\iota_\frakp:\calO \hookrightarrow \M_2(\Z_{F,\frakp})$.  The elements $\alpha_a$ in (\ref{deltaa}) are then generators with totally positive reduced norm of the left ideals 
\begin{equation} \label{Iadef}
I_a=\calO\iota_\frakp^{-1}\begin{pmatrix} x & y \\ 0 & 0 \end{pmatrix} + \calO\frakp
\end{equation}
and are obtained by principalizing the ideals $I_a$; here again we use strong approximation and the hypothesis that $F$ has strict class number $1$.  Then for each $a \in \PP^1(\F_p)$ and each $\gamma \in G$, we compute the permutation $\gamma^*$ \cite[Algorithm 5.8]{GV} and the element $\delta_a=\alpha_a\gamma\alpha_{\gamma^*a}^{-1} \in \Gamma(1)$ as in (\ref{deltaa}).  Using the solution to the word problem, we then write $\delta_a$ as a word in the generators $G$ for $\Gamma(1)$, and then for a basis of crossed homomorphisms $f$ we compute $f \slsh{} T_\frakp$ by computing $(f \slsh{} T_\frakp)(\gamma) \in V(\C)$ for each $\gamma \in G$ as in (\ref{E:heckeformula}).
In a similar way, we compute the action of complex conjugation and the Atkin-Lehner involutions.  We then decompose the space $H^1(\Gamma,V(\C))$ under the action of these operators into Hecke irreducible subspaces, and from this we compute the systems of Hecke eigenvalues using linear algebra.

\section{The indefinite method with arbitrary class number}

In this section, we show how to extend the method introduced in the previous section to the case where $F$ has arbitrary class number \cite[Remark 3.11]{GV}.  We refer the reader to Hida \cite{HidaHilMod} for further background.  

\subsection{Setup}

We carry over the notation from Section 1.  Recall that $\calO=\calO_0(\frakN)$ is an Eichler order of level $\frakN$ in the maximal order $\calO(1) \subset B$.

Let $\calH^{\pm} = \{z \in \C : \impart(z) \neq 0\} = \C \setminus \R$ be the union of the upper and lower half-planes.  Then via $\iota_\infty$, the group $B^\times$ acts on $\calH^{\pm}$ by linear fractional transformations.  

In this generality, we find it most elucidating to employ adelic notation.  Let $\Zhat = \varprojlim_{n} \Z/n\Z$ and let $\widehat{\phantom{x}}$ denote tensor with $\Zhat$ over $\Z$.  Consider the double coset
\[ X(\C)=B^\times \backslash (\calH^{\pm} \times \Bhat^\times / \calOhat^\times), \]
where $B^\times$ acts on $\Bhat^\times/\calOhat^\times$ by left multiplication via the diagonal embedding.  Then $X(\C)$ has the structure of a complex analytic space \cite{Deligne} which fails to be compact if and only if $B \cong \M_2(\Q)$, corresponding to the classical case of elliptic modular forms---higher class number issues do not arise in this case, so from now we assume that $B$ is a division ring.

We again write $S_k^B(\frakN)$ for the finite-dimensional $\C$-vector space of quaternionic modular forms of weight $k$ and level $\frakN$: here, again roughly speaking, a quaternionic modular form of weight $k \in (2\Z_{>0})^n$ and level $\frakN$ for $B$ is an analytic function 
\[ f: \calH^{\pm} \times \Bhat^\times / \calOhat^\times \to W_k(\C) \]
which is invariant under the weight $k$ action of $B^\times$, with $W_k(\C)$ as in Section 1.  

\subsection{Decomposing the double coset space}

By Eichler's theorem of norms, we have $\nrd(B^\times)=F_{(+)}^\times$ where 
\[ F_{(+)}^\times = \{ a \in F^\times : v_i(a)>0\text{ for $i=2,\dots,n$}\} \]
is the subgroup of elements of $F$ which are positive at all real places which are ramified in $B$.  In particular, $B^\times/B_+^\times \cong \Z/2\Z$, where 
\[ B_+^\times=\{\gamma \in B^\times : v_1(\nrd(\gamma))>0\}=\{ \gamma \in B : \nrd(\gamma) \in F_+^\times\}. \] 
The group $B_+^\times$ acts on the upper half-plane $\calH$, therefore we may identify
\[ X(\C) = B_+^\times \backslash (\calH \times \Bhat^\times / \calOhat^\times). \]
Now we have a natural (continuous) projection map
\[ X(\C) \to B_+^\times \backslash \Bhat^\times / \calOhat^\times, \]
and by strong approximation \cite[Theor\`eme III.4.3]{Vigneras} the reduced norm gives a bijection
\begin{equation} \label{strongapprox}
\nrd:B_+^\times \backslash \Bhat^\times / \calOhat^\times \xrightarrow{\sim} F_+^\times \backslash \Fhat^\times / \ZFhat^\times \cong \Cl^+ \Z_F,
\end{equation}
where $\Cl^+ \Z_F$ denotes the strict class group of $\Z_F$, i.e.\ the ray class group of $\Z_F$ with modulus equal to the product of all real (infinite) places of $F$.

The space $X(\C)$ is therefore the disjoint union of Riemann surfaces indexed by $\Cl^+ \Z_F$, which we identify explicitly as follows.  Let the ideals $\frakb \subset \Z_F$ form a set of representatives for $\Cl^+ \Z_F$, and let $\bhat \in \ZFhat$ be such that $\bhat\,\ZFhat \cap \Z_F = \frakb$.  For expositional simplicity, choose $\frakb=\Z_F$ and $\betahat=\widehat{1}$ for the representatives of the trivial class.  
By strong approximation (\ref{strongapprox}), there exists $\betahat \in \Bhat^\times$ such that $\nrd(\betahat)=\bhat$.  Therefore
\begin{equation} \label{breakup}
 X(\C) = \bigsqcup_{[\frakb]} B_+^\times (\calH \times \betahat \calOhat^\times).
\end{equation}
We have a map
\begin{align*} 
B_+^\times(\calH \times \betahat \calOhat^\times) &\to \calO_{\betahat,+}^\times \backslash \calH \\
(z,\betahat \calOhat^\times) &\mapsto z
\end{align*}
where $\calO_{\betahat}=\betahat \calOhat \betahat^{-1} \cap B$ and $\calO_{\betahat,+}^\times = \calO_{\betahat}^\times \cap B_+^\times$, so that $\calO_{\widehat{1}}=\calO$.

For each $\betahat$, let $\Gamma_{\betahat}=\iota_\infty\bigl(\calO_{\betahat,+}^\times/ \Z_F^\times\bigr) \subset \PGL_2(\R)^+$.  Then the Eichler-Shimura isomorphism on each component in (\ref{breakup}) gives an identification of Hecke modules
\begin{equation} \label{woahSk0}
S_k^B(\frakN) \xrightarrow{\sim} \bigoplus_{\betahat} H^1(\Gamma_{\betahat}, V_k(\C))^+,
\end{equation}
where ${}^+$ denotes the $+1$-eigenspace for complex conjugation.  For each $\betahat$, let $\calO(1)_{\betahat} = \betahat \calO(1) \betahat^{-1} \cap B$ be the maximal order containing the Eichler order $\calO_{\betahat}$, and let $\Gamma(1)_{\betahat} = \iota_\infty(\calO(1)_{\betahat,+}^\times/\Z_F^\times)$.  Further, let $V_{\betahat}(\C)=\Coind_{\Gamma_{\betahat}}^{\Gamma(1)_{\betahat}} V_k(\C)$.  Then Shapiro's lemma applied to each summand in (\ref{woahSk0}) gives
\begin{equation} \label{woahSk}
S_k^B(\frakN) \xrightarrow{\sim} \bigoplus_{\betahat} H^1(\Gamma(1)_{\betahat}, V_{\betahat}(\C))^+.
\end{equation}

\subsection{Hecke operators}

In the description (\ref{woahSk}), the Hecke operators $T_\frakp$ act on $\bigoplus_{\betahat} H^1(\Gamma(1)_{\betahat}, V_{\betahat}(\C))$ in the following way.  Let $\frakp$ be a prime ideal of $\Z_F$ with $\frakp \nmid \frakD\frakN$, and let $\phat \in \ZFhat$ be such that $\phat\,\ZFhat \cap \Z_F = \frakp$.  
% We suppose further that $\frakp$ is coprime to all ideals $\frakb$ chosen as representatives.  
We consider the $\betahat'$-summand in (\ref{woahSk}), corresponding to the ideal class $[\frakb']$.  Let $f:\Gamma(1)_{\betahat'} \to V_{\betahat'}(\C)$ be a crossed homomorphism: we will then obtain a new crossed homomorphism $f \slsh{} T_\frakp : \Gamma(1)_{\betahat} \to V_{\betahat}(\C)$, where $\betahat$ corresponds to the ideal class of $[\frakp\frakb']$ among the explicit choices made above.

Let $\varpihat \in \calOhat_{\betahat}$ be such that $\nrd(\varpihat)=\phat$.  Then there are elements $\gammahat_a \in \calOhat_{\betahat}$, indexed by $a \in \PP^1(\F_\frakp)$, such that
\begin{equation} \label{decompOstar}
\calOhat_{\betahat}^\times \varpihat \calOhat_{\betahat}^\times = \bigsqcup_{a \in \PP^1(\F_\frakp)} \calOhat_{\betahat}^\times \alphahat_a
\end{equation}
where $\alphahat_a = \varpihat \gammahat_a$.  

Let $\gamma \in \Gamma_{\betahat}$.  Extending (\ref{decompOstar}) to $\calOhat(1)_{\betahat}^\times$, we conclude that there exist unique elements $\deltahat_a \in \calOhat(1)_{\betahat}^\times$ and a unique permutation $\gamma^*$ of $\PP^1(\F_\frakp)$ such that
\[ \alphahat_a \gamma = \deltahat_a \alphahat_{\gamma^* a} \]
for $a \in \PP^1(\F_\frakp)$.  Thus we have
\[ (\betahat' \betahat^{-1} \alphahat_a) \gamma = (\betahat' \betahat^{-1}) \deltahat_a \alphahat_{\gamma^* a} =   \deltahat_a' (\betahat' \betahat^{-1} \alphahat_{\gamma^* a}). \]
where $\deltahat_a'=(\betahat' \betahat^{-1}) \deltahat_a (\betahat' \betahat^{-1})^{-1}$.  

Recall that $\betahat' \calOhat$ has left order $\calOhat_{\betahat'}$ and similarly $\calOhat \betahat^{-1}$ has right order $\calOhat_{\betahat}$.  Therefore, we may consider the left $\calOhat_{\betahat'}$-ideal
\begin{equation} \label{leftideal} 
\calOhat_{\betahat'} \betahat' \calOhat \betahat^{-1} \calOhat_{\betahat} \alphahat_a
\end{equation}
noting that the left and right orders in each case match up, so the product is compatible.  Next, 
recall that the elements $\betahat'$, $\betahat$, $\varpihat$ have reduced norms corresponding to the ideal classes $[\frakb']$, $[\frakp \frakb']$, and $[\frakp]$, respectively.  Thus the reduced norm of the left ideal (\ref{leftideal}) has a trivial ideal class.  Therefore, by strong approximation (applied now to left ideals of the order $\calO_{\betahat'}$), for each $a \in \PP^1(\F_\frakp)$, there exist elements $\pi_a' \in \calO_{\betahat'} \cap B_+^\times$ such that
\[ \calOhat_{\betahat'} \betahat' \betahat^{-1} \alphahat_a \cap B = \calO_{\betahat'} \pi_a'. \]
Hence there exists a unique permutation $\gamma^*$ of $\PP^1(\F_\frakp)$ such that
\[ \pi_a' \gamma = \delta_a' \pi_{\gamma^* a}' \]
with $\delta_a \in \calO_{\betahat',+}^\times$.  The new crossed homomorphism $f \slsh{} T_\frakp : \Gamma_{\betahat} \to V_{\betahat}(\C)$ is then defined by the formula
\[ (f \slsh{} T_\frakp)(\gamma) = \sum_{a\in\PP^1(\F_\frakp)}f(\delta_a')^{\pi_a'} \]
for $\gamma \in \Gamma_{\betahat}$.

\subsection{Complex conjugation and Atkin-Lehner involutions}

We now define an operator $W_\infty$ which acts by complex conjugation.  Let $\Cl^{(+)} \Z_F$ denote the ray class group of $\Z_F$ with modulus equal to the real (infinite) places of $F$ which are ramified in $B$.  Then we have a natural map
$\Cl^+ \Z_F \to \Cl^{(+)} \Z_F$; this map is an isomorphism if and only if there exists a unit $u \in \Z_F^\times$ which satisfies $v_1(u)<0$ and $v_i(u)>0$ for the other real places $v_i$ ($i=2,\dots,n$) of $F$, otherwise the kernel of this map is isomorphic to $\Z/2\Z$.  Let $[\frakm] \in \Cl^+ \Z_F$ generate the kernel of this map.  

Let $f:\Gamma(1)_{\betahat'} \to V_{\betahat'}$ be a crossed homomorphism, and let $\betahat$ correspond to the ideal class $[\frakb' \frakm^{-1}]$; we will define the complex conjugate crossed homomorphism $(f \slsh{} W_\infty) : \Gamma(1)_{\betahat} \to V_{\betahat}(\C)$.  The left $\calO_{\betahat'}$-ideal $\calOhat_{\betahat'} \betahat' \calOhat \betahat^{-1} \cap B$ has reduced norm corresponding to the ideal class $[\frakm] \in \Cl^+ \Z_F$, so there exists a generator $\mu' \in \calO_{\betahat'}$ of this ideal such that $v_1(\nrd(\mu'))<0$ but $v_i(\nrd(\mu'))>0$ for $i=2,\dots,n$.  Then given $\gamma \in \Gamma(1)_{\betahat}$, we define
\[ (f \slsh{} W_\infty)(\gamma) = f(\mu' \gamma \mu'^{-1})^{\mu'}. \]

Finally, we define the Atkin-Lehner involutions $W_{\frakp^e}$ for $\frakp^e \parallel \frakD\frakN$.  Let $\frakp$ correspond to $\phat \in \ZFhat$.  Then there exists an element $\pihat \in \calO_{\betahat}$ which generates the unique two-sided ideal of $\calO_{\betahat}$ of reduced norm generated by $\phat^e$.  The element $\pihat$ normalizes $\calO_{\betahat}$ and $\pihat^2 \in \calO_{\betahat}^\times \Fhat^\times$.  Let $\betahat$ correspond to the ideal class $[\frakp \frakb']$.  Then as above, by strong approximation there exists an element $\mu' \in \calO_{\betahat'} \cap B_+^\times$ such that $\calO_{\betahat'} \betahat' \betahat \pihat \cap B = \calO_{\betahat'} \mu'$.  
Given $f:\Gamma(1)_{\betahat'} \to V_{\betahat'}$, we then define $(f \slsh{} W_{\frakp^e}) : \Gamma(1)_{\betahat} \to V_{\betahat}(\C)$ by
\[ (f \slsh{} W_\frakp)(\gamma) = f(\mu' \gamma \mu'^{-1})^{\mu'} \]
for $\gamma \in \Gamma(1)_{\betahat}$.

\section{Algorithmic methods}

In this section, we take the adelic description of Section 2 and show how to compute with it explicitly, proving Theorem 1.  

Our algorithm takes as input a totally real field $F$ of degree $[F:\Q]=n$, a quaternion algebra $B$ over $F$ split at a unique real place, an ideal $\frakN \subset \Z_F$ coprime to the discriminant $\frakD$ of $B$, a vector $k \in (2\Z_{>0})^n$, and a prime $\frakp \nmid \frakD\frakN$, and outputs the matrix of the Hecke operator $T_\frakp$ acting on the space $H= \bigoplus_{\betahat} H^1\bigl(\Gamma(1)_{\betahat}, V_{\betahat}(\C)\bigr)^+$ (in the notation of Section 2) with respect to some fixed basis which does not depend on $\frakp$.  From these matrices, one decomposes the space $H$ into Hecke-irreducible subspaces by the techniques of basic linear algebra.

Our algorithm follows the form given in the overview in Section 1, so we describe our algorithm in steps, with a description of each step along the way.

\medskip\noindent \textbf{Step 1} (\textsf{Compute a splitting field}):  Let $K \hookrightarrow \C$ be a Galois number field containing $F$ which splits $B$: for example, we can take the normal closure of any quadratic field contained in $B$.  Since all computations then occur inside $K \subset \C$, we may work then with coefficient modules over $K$ using exact arithmetic.  (This step is only necessary if $k$ is not parallel weight $2$, for otherwise the action of $B^\times$ factors through $K=\Q$.)

% \[ V_{\betahat}(K) = \Coind_{\Gamma_{\betahat}}^{\Gamma(1)_{\betahat}} V_k(K) \]
% where $V_k(K)$ is as in (\ref{Vk}) but with coefficients in $K$, since $V_{\betahat}(\C) = V_{\betahat}(K) \otimes_K \C$.  The $K$-vector spaces $V_{\betahat}(K)$ are then equipped with an action of $B^\times$ which can be represented using exact arithmetic over $K$.  

\medskip\noindent \textbf{Step 2} (\textsf{Compute ideal class representatives}):  Compute a set of representatives $[\frakb]$ for the strict class group $\Cl^+ \Z_F$ with each $\frakb$ coprime to $\frakp\frakD\frakN$.  (See Remark \ref{bcoprime} below.)

Compute a maximal order $\calO(1) \subset B$.  For each representative ideal $\frakb$, compute a right $\calO(1)$-ideal $J_{\frakb}$ such that $\nrd(J_{\frakb}) = \frakb$ and let $\calO(1)_\frakb$ be the left order of $J_{\frakb}$.  (In the notation of Section 2, the right $\calO(1)$-ideals $J_{\frakb}$ represent the elements $\betahat$, and $\calO(1)_\frakb=\calO(1)_{\betahat}$.)

\medskip\noindent \textbf{Step 3} (\textsf{Compute presentations for the unit groups}):  Compute an embedding $\iota_\infty: B \hookrightarrow \M_2(\R)$ corresponding to the unique split real place.  

For each $\frakb$, compute a finite presentation for $\Gamma(1)_\frakb = \iota_{\infty}(\calO(1)_{\frakb,+}^\times/\Z_F^\times)$ consisting of a (minimal) set of generators $G_\frakb$ and relations $R_\frakb$ together with a solution to the word problem for the computed presentation \cite{V-fd}.  (Note that the algorithm stated therein \cite[Theorem 3.2]{V-fd} is easily extended from units of reduced norm $1$ to totally positive units.)

For efficiency, we start by computing such a presentation with generators $G$ associated to the order $\calO(1)$ and then for each order $\calO(1)_{\frakb}$ we begin with the elements in hand formed by short products of elements in $G$ which happen to lie in $\calO(1)_{\frakb}$ (to aid in the search for units \cite[Algorithm 3.2]{V-fd}; note that $\calO(1) \cap \calO(1)_{\frakb}$ is an Eichler order of level $\frakb$ in $\calO(1)_{\frakb}$).

\medskip\noindent \textbf{Step 4} (\textsf{Compute splitting data}): Compute a splitting 
\[ \iota_{\frakN} : \calO(1) \hookrightarrow \calO(1) \otimes_{\Z_F} \Z_{F,\frakN} \cong \M_2(\Z_{F,\frakN}). \]  
Note that since $\frakb$ is coprime to $\frakN$, we have $\calO(1) \otimes \Z_{F,\frakN} = \calO(1)_\frakb \otimes \Z_{F,\frakN}$ for all $\frakb$, so $\iota_{\frakN}$ also gives rise to a splitting for each $\calO(1)_\frakb$.  For each $\frakb$, compute the Eichler order $\calO_\frakb \subset \calO(1)_\frakb$ of level $\frakN$ with respect to $\iota_\frakN$.  

Next, for each $\frakb$, compute representatives for the left cosets of the group $\Gamma_{\frakb} = \iota_\infty(\calO_{\frakb,+}^\times/\Z_F^\times)$ inside $\Gamma(1)_{\frakb}$ \cite[Algorithm 6.1]{GV}.  Finally, identify 
\[ V(K)_{\frakb}= \Coind_{\Gamma_\frakb}^{\Gamma(1)_\frakb} V_k(K) \] 
as a $K$-vector space given by copies of $V_k(K)$ indexed by these cosets, and compute the permutation action of the representatives of these cosets on this space.

In practice, it is more efficient to identify the above coset representatives with elements of $\PP^1(\Z_F/\N)$ and thereby work directly with the coefficient module $V(K)_{\frakb} \cong K[\PP^1(\Z_F/\frakN)] \otimes V_k(K)$.

\medskip\noindent \textbf{Step 5} (\textsf{Compute a basis for cohomology}):  Identify the space of crossed homomorphisms $\bigoplus_{\frakb} Z^1(\Gamma(1)_\frakb, V(K)_\frakb)$ with its image under the inclusion
\begin{align*} 
Z^1(\Gamma_\frakb,V(K)_\frakb) &\to \bigoplus_{g \in G_\frakb} V(K)_\frakb \\
f &\mapsto (f(g))_{g\in G_\frakb}
\end{align*}
consisting of those $f \in \bigoplus_{g \in G_\frakb} V(K)_\frakb$ which satisfy the  relations $f(r)=0$ for $r \in R_\frakb$.  Compute the space of principal crossed homomorphisms $B^1(\Gamma(1)_\frakb, V(K)_\frakb)$ in a similar way, and thereby compute using linear algebra a $K$-basis for the quotient $H^1(\Gamma(1)_\frakb, V(K)_\frakb) = Z^1(\Gamma(1)_\frakb, V(K)_\frakb)/B^1(\Gamma(1)_\frakb, V(K)_\frakb)$ for each $\frakb$.

Let $H=\bigoplus_{\frakb} H^1(\Gamma(1)_\frakb, V(K)_\frakb)$.  

\medskip\noindent \textbf{Step 6} (\textsf{Compute representatives for left ideal classes}):  Compute a splitting $\iota_\frakp: \calO(1) \hookrightarrow \M_2(\Z_{F,\frakp})$.  For each ideal $\frakb'$, perform the following steps.  

First, compute the ideal $\frakb$ with ideal class $[\frakb]=[\frakp \frakb']$.  Compute the left ideals 
\[ I_a =\calO\iota_\frakp^{-1}\begin{pmatrix} x & y \\ 0 & 0 \end{pmatrix} + \calO\frakp \]
indexed by the elements $a=(x:y) \in \PP^1(\F_\frakp)$ and then compute the left $\calO_{\frakb'}$-ideals $I_a'=J_{\frakb'} \overline{J}_\frakb I_a$.  

Compute totally positive generators $\pi_a' \in \calO_{\frakb'} \cap B_+^\times$ for $\calO_{\frakb'} \pi_a' = I_a'$ \cite{KV}.

Now, for each $\gamma \in G_\frakb$, compute the permutation $\gamma^*$ of $\PP^1(\F_\frakp)$ \cite[Algorithm 5.8]{GV} and then the elements $\delta_a'=\pi_a'\gamma\pi_{\gamma^*a}'^{-1}$ for $a \in \PP^1(\F_\frakp)$; write each such element $\delta_a'$ as a word in $G_\frakb'$ and from the formula
\[ (f \slsh{} T_\frakp)(\gamma) = \sum_{a\in\PP^1(\F_\frakp)}f(\delta_a')^{\pi_a'} \]
with $f$ in a basis for the $\frakb'$-component of cohomology as in Step 5 compute the induced crossed homomorphism $f \slsh{} T_\frakp$ in the $\frakb$-component.  

\medskip\noindent \textbf{Step 7} (\textsf{Compute the blocks of the intermediate matrix}): Assemble the matrix $T$ with rows and columns indexed as in Step 5 with blocks in the $(\frakb,\frakb')$ position given by the output of Step 6: this matrix describes the action of $T_\frakp$ on $H$.

\medskip\noindent \textbf{Step 8} (\textsf{Decompose $H$ into $\pm$-eigenspaces for complex conjugation}): Determine the representative ideal $\frakm$ (among the ideals $\frakb$) which generates the kernel of the map $\Cl^+ \Z_F \to \Cl^{(+)} \Z_F$.  

For each ideal $\frakb'$, perform the following steps.  Compute the ideal $\frakb$ such that $[\frakb]=[\frakb'\frakm^{-1}]$, and compute a generator $\mu'$ with $\calO_{\frakb'} \mu' = J_{\frakb'} \overline{J}_{\frakb}$ such that $v(\nrd(\mu'))<0$.  
For each $\gamma \in G_\frakb$, from the formula
\[ (f \slsh{} W_\infty)(\gamma) = f(\mu'\gamma \mu'^{-1})^{\mu'}, \]
for $f$ in a basis for the $\frakb'$-component of cohomology as in Step 5 compute the induced crossed homomorphism $f \slsh{} T_\frakp$ in the $\frakb$-component.  

Assemble the matrix with blocks in the $(\frakb,\frakb')$ position given by this output: this matrix describes the action of complex conjugation $W_\infty$ on $H$.  Compute a $K$-basis for the $+1$-eigenspace $H^+$ of $H$ for $W_\infty$.  Finally, compute the matrix $T^+$ giving the action of $T_\frakp$ restricted to $H^+$ and return $T^+$.

\medskip This completes the description of the algorithm.  

In a similar way, one computes the Atkin-Lehner involutions, replacing Step 6 with the description given in Section 2.4, similar to the computation of complex conjugation in Step 8.

\begin{remark}
Note that Steps 1 through 3 do not depend on the prime $\frakp$ nor the level $\frakN$ and Steps 4, 5, and 8 do not depend on the prime $\frakp$, so these may be precomputed for use in tabulation.
\end{remark}

\begin{remark} \label{bcoprime}
To arrange uniformly that the ideals $\frakb$ representing the classes in $\Cl^+ \Z_F$ are coprime to the prime $\frakp$ in advance for many primes $\frakp$, one has several options.  One possibility is to choose suitable ideals $\frakb$ of large norm in advance.  Another option is to make suitable modifications ``on the fly'': if $\frakp$ is not coprime to $\frakb$, we simply choose a different ideal $\frakc$ coprime to $\frakp$ with $[\frakb]=[\frakc]$, a new ideal $J_\frakc$ with $\nrd(J_\frakc)=\frakc$, and compute an element $\nu \in \calO_\frakb$ such that $\nu \calO_\frakb \nu^{-1} = \calO_\frakc$.  Conjugating by $\nu$ where necessary, one can then transport the computations from one order to the other so no additional computations need to take place.
\end{remark}

\section{Examples}

In this section, we compute with two examples to demonstrate the algorithm outlined in Section 3.  Throughout, we use the computer system \textsf{Magma} \cite{Magma}.

Our first and most detailed example is concerned with the smallest totally real cubic field $F$ with the property that the dimension of the space of Hilbert cusp forms of parallel weight $2$ and level $(1)$ is greater than zero and the strict class number of $F$ is equal to $2$.  This field is given by $F=\Q(w)$ where $w$ satisfies the equation $f(w)=w^3-11w-11=0$.  The discriminant of $F$ is equal to $2057=11^2 17$, and $\Z_F=\Z[w]$.  The roots of $f$ in $\R$ are $-2.602\dots$, $-1.131\dots$, and $3.73\dots$, and we label the real places $v_1,v_2,v_3$ of $F$ into $\R$ according to this ordering.

We define the \emph{sign} of $a \in F$ to be the triple $\sgn(a)=(\sgn(v_i(a)))_{i=1}^{3} \in \{\pm 1\}^3$.  
The unit group of $F$ is generated by the elements $-1$, $w+1$ with $\sgn(w+1)=(1,-1,-1)$, and the totally positive unit $-w^2+2w+12$.

We begin by finding a quaternion algebra $B$ with $\frakD=\Z_F$ which is ramified at all but one real place \cite[Algorithm 4.1]{GV}.  We find the algebra $B=\quat{w+1,-1}{F}$ ramified only at $v_1$ and $v_2$, generated by $i,j$ subject to $i^2=w+1$, $j^2=-1$, and $ji=-ij$.

For forms of parallel weight $2$, \textbf{Step 1} is trivial: we can take $K=\Q$.  

Next, in \textbf{Step 2} we compute ideal class representatives.  The nontrivial class in $\Cl^+(\Z_F)$ is represented by the ideal $\frakb=(w^2-2w-6)\Z_F$, which is principal but does not possess a totally positive generator, since $\sgn(-w^2+2w+6)=(-1,1,-1)$ and there is no unit of $\Z_F$ with this sign.  We note that $\N(\frakb)=7$.

Next, we compute a maximal order $\calO=\calO(1)$; it is generated over $\Z_F$ by $i$ and the element $k=(1+(w^2+1)i+ij)/2$.  Next, we find that the right $\calO$-ideal $J_\frakb$ generated by $w^2-2w-6$ and the element $(5+(w^2+5)i+ij)/2=2+2i+k$ has $\nrd(J_\frakb)=\frakb$.

Next, in \textbf{Step 3} we compute presentations for the unit groups.  We take the splitting 
\begin{align*}
B &\hookrightarrow \M_2(\R) \\
i,j &\mapsto \begin{pmatrix} s & 0 \\ 0 & -s \end{pmatrix}, \begin{pmatrix} 0 & 1 \\ -1 & 0 \end{pmatrix} 
\end{align*}
where $s=\sqrt{v_3(w+1)}$.  We then compute a fundamental domain for $\Gamma=\Gamma(1)$ \cite{V-fd}, given below.

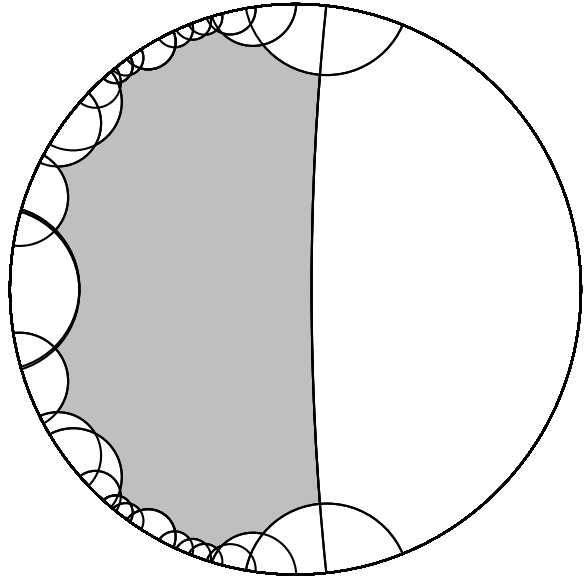
\begin{center}
\psset{unit=1.5in}
\begin{pspicture}(-1,-1)(1,1)
\pscircle[fillstyle=solid,fillcolor=lightgray](0,0){1}

\psclip{\pscircle(0,0){1}} \pscircle[fillstyle=solid,fillcolor=white](9.52632,-0.000000){9.47368} \endpsclip
\psclip{\pscircle(0,0){1}} \pscircle[fillstyle=solid,fillcolor=white](9.52632,-0.000000){9.47368} \endpsclip
\psclip{\pscircle(0,0){1}} \pscircle[fillstyle=solid,fillcolor=white](0.109426,1.03580){0.291287} \endpsclip
\psclip{\pscircle(0,0){1}} \pscircle[fillstyle=solid,fillcolor=white](-0.148315,1.00097){0.154706} \endpsclip
\psclip{\pscircle(0,0){1}} \pscircle[fillstyle=solid,fillcolor=white](-0.227238,0.978175){0.0920039} \endpsclip
\psclip{\pscircle(0,0){1}} \pscircle[fillstyle=solid,fillcolor=white](-0.316782,0.950756){0.0654914} \endpsclip
\psclip{\pscircle(0,0){1}} \pscircle[fillstyle=solid,fillcolor=white](-0.358943,0.935825){0.0678953} \endpsclip
\psclip{\pscircle(0,0){1}} \pscircle[fillstyle=solid,fillcolor=white](-0.421475,0.909378){0.0678953} \endpsclip
\psclip{\pscircle(0,0){1}} \pscircle[fillstyle=solid,fillcolor=white](-0.513834,0.863631){0.0994178} \endpsclip
\psclip{\pscircle(0,0){1}} \pscircle[fillstyle=solid,fillcolor=white](-0.513834,0.863631){0.0994178} \endpsclip
\psclip{\pscircle(0,0){1}} \pscircle[fillstyle=solid,fillcolor=white](-0.592241,0.808417){0.0654914} \endpsclip
\psclip{\pscircle(0,0){1}} \pscircle[fillstyle=solid,fillcolor=white](-0.628515,0.780455){0.0643415} \endpsclip
\psclip{\pscircle(0,0){1}} \pscircle[fillstyle=solid,fillcolor=white](-0.628515,0.780455){0.0643415} \endpsclip
\psclip{\pscircle(0,0){1}} \pscircle[fillstyle=solid,fillcolor=white](-0.697694,0.722280){0.0920039} \endpsclip
\psclip{\pscircle(0,0){1}} \pscircle[fillstyle=solid,fillcolor=white](-0.775487,0.654707){0.173267} \endpsclip
\psclip{\pscircle(0,0){1}} \pscircle[fillstyle=solid,fillcolor=white](-0.829191,0.579979){0.154706} \endpsclip
\psclip{\pscircle(0,0){1}} \pscircle[fillstyle=solid,fillcolor=white](-0.962779,0.321059){0.173267} \endpsclip
\psclip{\pscircle(0,0){1}} \pscircle[fillstyle=solid,fillcolor=white](-1.04155,0.00447809){0.291287} \endpsclip
\psclip{\pscircle(0,0){1}} \pscircle[fillstyle=solid,fillcolor=white](-1.04155,-0.00447809){0.291287} \endpsclip
\psclip{\pscircle(0,0){1}} \pscircle[fillstyle=solid,fillcolor=white](-0.962779,-0.321059){0.173267} \endpsclip
\psclip{\pscircle(0,0){1}} \pscircle[fillstyle=solid,fillcolor=white](-0.829191,-0.579979){0.154706} \endpsclip
\psclip{\pscircle(0,0){1}} \pscircle[fillstyle=solid,fillcolor=white](-0.775487,-0.654707){0.173267} \endpsclip
\psclip{\pscircle(0,0){1}} \pscircle[fillstyle=solid,fillcolor=white](-0.697694,-0.722280){0.0920039} \endpsclip
\psclip{\pscircle(0,0){1}} \pscircle[fillstyle=solid,fillcolor=white](-0.628515,-0.780455){0.0643415} \endpsclip
\psclip{\pscircle(0,0){1}} \pscircle[fillstyle=solid,fillcolor=white](-0.628515,-0.780455){0.0643415} \endpsclip
\psclip{\pscircle(0,0){1}} \pscircle[fillstyle=solid,fillcolor=white](-0.592241,-0.808417){0.0654914} \endpsclip
\psclip{\pscircle(0,0){1}} \pscircle[fillstyle=solid,fillcolor=white](-0.513834,-0.863631){0.0994178} \endpsclip
\psclip{\pscircle(0,0){1}} \pscircle[fillstyle=solid,fillcolor=white](-0.513834,-0.863631){0.0994178} \endpsclip
\psclip{\pscircle(0,0){1}} \pscircle[fillstyle=solid,fillcolor=white](-0.421475,-0.909378){0.0678953} \endpsclip
\psclip{\pscircle(0,0){1}} \pscircle[fillstyle=solid,fillcolor=white](-0.358943,-0.935825){0.0678953} \endpsclip
\psclip{\pscircle(0,0){1}} \pscircle[fillstyle=solid,fillcolor=white](-0.316782,-0.950756){0.0654914} \endpsclip
\psclip{\pscircle(0,0){1}} \pscircle[fillstyle=solid,fillcolor=white](-0.227238,-0.978175){0.0920039} \endpsclip
\psclip{\pscircle(0,0){1}} \pscircle[fillstyle=solid,fillcolor=white](-0.148315,-1.00097){0.154706} \endpsclip
\psclip{\pscircle(0,0){1}} \pscircle[fillstyle=solid,fillcolor=white](0.109426,-1.03580){0.291287} \endpsclip

\psclip{\pscircle(0,0){1}} \pscircle(9.52632,-0.000000){9.47368} \endpsclip
\psclip{\pscircle(0,0){1}} \pscircle(9.52632,-0.000000){9.47368} \endpsclip
\psclip{\pscircle(0,0){1}} \pscircle(0.109426,1.03580){0.291287} \endpsclip
\psclip{\pscircle(0,0){1}} \pscircle(-0.148315,1.00097){0.154706} \endpsclip
\psclip{\pscircle(0,0){1}} \pscircle(-0.227238,0.978175){0.0920039} \endpsclip
\psclip{\pscircle(0,0){1}} \pscircle(-0.316782,0.950756){0.0654914} \endpsclip
\psclip{\pscircle(0,0){1}} \pscircle(-0.358943,0.935825){0.0678953} \endpsclip
\psclip{\pscircle(0,0){1}} \pscircle(-0.421475,0.909378){0.0678953} \endpsclip
\psclip{\pscircle(0,0){1}} \pscircle(-0.513834,0.863631){0.0994178} \endpsclip
\psclip{\pscircle(0,0){1}} \pscircle(-0.513834,0.863631){0.0994178} \endpsclip
\psclip{\pscircle(0,0){1}} \pscircle(-0.592241,0.808417){0.0654914} \endpsclip
\psclip{\pscircle(0,0){1}} \pscircle(-0.628515,0.780455){0.0643415} \endpsclip
\psclip{\pscircle(0,0){1}} \pscircle(-0.628515,0.780455){0.0643415} \endpsclip
\psclip{\pscircle(0,0){1}} \pscircle(-0.697694,0.722280){0.0920039} \endpsclip
\psclip{\pscircle(0,0){1}} \pscircle(-0.775487,0.654707){0.173267} \endpsclip
\psclip{\pscircle(0,0){1}} \pscircle(-0.829191,0.579979){0.154706} \endpsclip
\psclip{\pscircle(0,0){1}} \pscircle(-0.962779,0.321059){0.173267} \endpsclip
\psclip{\pscircle(0,0){1}} \pscircle(-1.04155,0.00447809){0.291287} \endpsclip
\psclip{\pscircle(0,0){1}} \pscircle(-1.04155,-0.00447809){0.291287} \endpsclip
\psclip{\pscircle(0,0){1}} \pscircle(-0.962779,-0.321059){0.173267} \endpsclip
\psclip{\pscircle(0,0){1}} \pscircle(-0.829191,-0.579979){0.154706} \endpsclip
\psclip{\pscircle(0,0){1}} \pscircle(-0.775487,-0.654707){0.173267} \endpsclip
\psclip{\pscircle(0,0){1}} \pscircle(-0.697694,-0.722280){0.0920039} \endpsclip
\psclip{\pscircle(0,0){1}} \pscircle(-0.628515,-0.780455){0.0643415} \endpsclip
\psclip{\pscircle(0,0){1}} \pscircle(-0.628515,-0.780455){0.0643415} \endpsclip
\psclip{\pscircle(0,0){1}} \pscircle(-0.592241,-0.808417){0.0654914} \endpsclip
\psclip{\pscircle(0,0){1}} \pscircle(-0.513834,-0.863631){0.0994178} \endpsclip
\psclip{\pscircle(0,0){1}} \pscircle(-0.513834,-0.863631){0.0994178} \endpsclip
\psclip{\pscircle(0,0){1}} \pscircle(-0.421475,-0.909378){0.0678953} \endpsclip
\psclip{\pscircle(0,0){1}} \pscircle(-0.358943,-0.935825){0.0678953} \endpsclip
\psclip{\pscircle(0,0){1}} \pscircle(-0.316782,-0.950756){0.0654914} \endpsclip
\psclip{\pscircle(0,0){1}} \pscircle(-0.227238,-0.978175){0.0920039} \endpsclip
\psclip{\pscircle(0,0){1}} \pscircle(-0.148315,-1.00097){0.154706} \endpsclip
\psclip{\pscircle(0,0){1}} \pscircle(0.109426,-1.03580){0.291287} \endpsclip
\pscircle(0,0){1}
\end{pspicture}
\end{center}

We find that $\Gamma=\Gamma(1)$ is the free group on the generators $\alpha,\beta,\gamma_1,\dots,\gamma_7$ subject to the relations
\[ \gamma_1^2=\gamma_2^2=\gamma_3^3=\gamma_4^2=\gamma_5^3=\gamma_6^2=\gamma_7^2=
\alpha\beta\alpha^{-1}\beta^{-1}\gamma_1\cdots \gamma_7 = 1. \]
For example, we have
\[ 2\alpha=(w^2 - 14) + (2w^2 - 4w - 13)i + (-2w^2 + 5w + 9)j + (-4w^2 + 8w + 26)ij. \]
% and
% \begin{align*}
% 2\beta &=(-23w^2 - 86w - 68) + (4w^2 + 15w + 12)i + \\
% &\qquad (46w^2 + 172w + 136)j + (-23w^2 - 86w - 68)ij.
% \end{align*}

The groups $\Gamma$ and $\Gamma_\frakb$ have isomorphic presentations.  In particular, we note that both $\Gamma$ and $\Gamma_\frakb$ have genus $1$, so we conclude that $\dim S_2(1)=1+1=2$.  

We illustrate the computation of Hecke operators with the primes $\frakp_3=(w+2)\Z_F$ of norm $3$ and $\frakp_5=(w+3)\Z_F$ of norm $5$.  Note that $\frakp_3$ is nontrivial in $\Cl^+(\Z_F)$ whereas $\frakp_5$ is trivial.  

Step \textbf{Step 4} requires no work, since we work with forms of level $(1)$.  In \textbf{Step 5} we compute with a basis for cohomology, and here we see directly that
\[ H^1(\Gamma,\Q) \cong \Hom(\Gamma,\Q) \cong \Z f_\alpha \oplus \Z f_\beta \]
where $f_\alpha, f_\beta$ are the characteristic functions for $\alpha$ and $\beta$.  We have a similar description for $H^1(\Gamma_\frakb,\Q)$. 

Next, in \textbf{Step 6} we compute representatives of the left ideal classes.  For $\frakp_3$, for example, for $I_{[1:0]} \subset \calO$ we find that $J_\frakb I_{[1:0]} = \calO_\frakb((w+1)+i+ij)$ and for $I_{[1:1]} \subset \calO_\frakb$ we have $\overline{J_\frakb} I_{[1:1]}=\calO(w+1-i+ij)$; we thereby find elements $\pi_a, \pi_a'$ for $a \in \PP^1(\F_{\frakp_3})$.  For the generators $\gamma=\alpha,\beta$ of $\calO$ and $\calO_\frakb$, we compute the permutations $\gamma^*$ of $\PP^1(\F_{\frakp_3})$; we find for example that $\alpha^*$ is the identity and
\[ \pi_{[1:0]}' \alpha = \delta_{[1:0]}' \pi_{[1:0]}' \]
with $\delta_{[1:0]}' \in \calO_\frakb$, namely,
\begin{align*}
14\delta_{[1:0]}' &= (7w^2 - 98) + (-23w^2 + 40w + 167)i + \\
&\qquad (-25w^2 + 59w + 103)j + (-2w^2 + 5w + 20)ij.
\end{align*}
We then write $\delta_{[1:0]}'$ as a word in the generators for $\Gamma_\frakb'$ of length $23$.  Repeating these steps (reducing a total of $64$ units), we assemble the block matrix in \textbf{Step 7} as the matrix
\[ T_{\frakp_3} \slsh{} H=
\begin{pmatrix}
0 & 0 & 2 & 0 \\
0 & 0 & 0 & 2 \\
2 & 0 & 0 & 0 \\
0 & 2 & 0 & 0
\end{pmatrix}. \]
In a similar way, we find that $T_{\frakp_5}$ is the identity matrix.

Finally, in \textbf{Step 8} we compute the action of complex conjugation.  Here we have simply $\mu=i$ (whereas $\mu_\frakb$ is more complicated), and thereby compute that
\[ W_\infty \slsh{} H =
\begin{pmatrix}
1 & 1 & 0 & 0 \\
0 & -1 & 0 & 0 \\
0 & 0 & 1 & 1 \\
0 & 0 & 0 & -1
\end{pmatrix}. \]
We verify that $W_\infty$ commutes with $T_{\frakp_3}$ (and $T_{\frakp_5}$).  We conclude that $T_{\frakp_3} \slsh{} H^+=\begin{pmatrix} 0 & 2 \\ 2 & 0 \end{pmatrix}$ and $T_{\frakp_5} \slsh{} H^+=\begin{pmatrix} 1 & 0 \\ 0 & 1 \end{pmatrix}$.

We then diagonalize the space $H^+$, which breaks up into two one-dimensional eigenforms $f$ and $g$, and compute several more Hecke operators: we list in Table 1 below a generator for the prime $\frakp$, its norm $\N\frakp$, and the Hecke eigenvalues $a_\frakp(f)$ and $a_\frakp(g)$ for the cusp forms $f,g$.  

\begin{table}[h]
\begin{equation} \label{tableF5l1} \notag
\setlength\arraycolsep{0.8em}
\begin{array}{c|c||cc}
\frakp & \N\!\frakp & a_\frakp(f) & a_\frakp(g) \\
\hline
w+2 & 3 & 2 & -2 \\
w+3 & 5 & 1 & 1 \\
2 & 8 & -5 & -5 \\
2w+7 & 9 & -2 & 2 \\
w & 11 & 0 & 0 \\
w^2-w-8 & 17 & -5 & 5 \\
w-3 & 17 & -5 & -5 \\
2w^2-5w-10 & 23 & 2 & -2 \\
w^2-3w-2 & 25 & -9 & -9 \\
w^2-6 & 29 & 9 & -9 \\
w+4 & 31 & -2 & -2 \\
2w^2-3w-16 & 37 & -3 & 3 \\
w^2-2w-9 & 41 & -5 & 5 \\
w^2+w-3 & 49 & -10 & 10 
\end{array}
\end{equation}
\textbf{Table 1}: Hecke eigenvalues for the Hilbert cusp forms for $F=\Q(w)$ with $w^3-11w-11=0$ of level $(1)$ and parallel weight $2$
\end{table}

We note that the primes generated by $w$ and $w-3$ are ramified in $F$.  

By work of Deligne \cite{Deligne}, the curves $X=X(1)$ and $X_\frakb$ are defined over the strict class field $F^+$ of $F$, and $\Gal(F^+/F)$ permutes them.  We compute that $F^+=F(\sqrt{-3w^2+8w+12})$.  Therefore the Jacobian $J_f$, corresponding to the cusp form $f$, is a modular elliptic curve over $F^+$ with $\#J(\F_\frakp)=\N\frakp+1-a_f(\frakp)$ with everywhere good reduction.  The form $g$ is visibly a quadratic twist of $f$ by the character corresponding to the extension $F^+/F$.  

Unfortunately, this curve does not have any apparent natural torsion structure which would easily allow for its identification as an explicit curve given by a sequence of coefficients \cite[\S 4]{DD}.

As a second and final example, we compute with a quaternion algebra defined over a quadratic field and therefore ramified at a finite prime.  We take $F=\Q(\sqrt{65})$, with $\Z_F=\Z[(1+\sqrt{65})/2]$.  The field $F$ has $\#\Cl(F)=\#\Cl^+(F)=2$.  We compute the space $S=S_2(\frakp_5)^{\textup{$\frakp_5$-new}}$ of Hilbert cuspidal new forms of parallel weight $2$ and level $\frakp_5$, where $\frakp_5$ is the unique prime in $\Z_F$ of norm $5$.  

We compute that $\dim S=10$, and that the space $S$ decomposes into Hecke-irreducible subspaces of dimensions $2,2,3,3$.  For example, the characteristic polynomial of $T_{\frakp_2}$ for $\frakp_2$ either prime above $2$ factors as
\[ (T^2-2T-1)(T^2+2T-1)(T^6 + 11T^4 + 31T^2 + 9). \]

\begin{remark} \label{dembeleoverlap}
By the Jacquet-Langlands correspondence, the space $S_2(\frakp_5)^{\textup{$\frakp_5$-new}}$ also occurs in the space of quaternionic modular forms for an Eichler order of level $\frakp_5$ in the definite quaternion algebra ramified at the the two real places of $F$ and no finite place, and therefore is amenable to calculation by the work of Demb\'el\'e and Donnelly.  We use this overlap to duplicate their computations (as well as ours) and thereby give some compelling evidence that the results are correct since they are computed in entirely different ways.
\end{remark}

\end{document}